\documentclass[12pt]{amsart}
\usepackage[totalwidth=480pt,totalheight=680pt]{geometry}
\input epsf.tex
\newtheorem{thm}{Theorem}[section]
\newtheorem{thm*}{Theorem}
\newtheorem{cor}[thm]{Corollary}
\newtheorem{lem}[thm]{Lemma}

\theoremstyle{definition}

\theoremstyle{remark}
\newtheorem{rem}[thm]{Remark}

\numberwithin{equation}{section}

\newcommand{\Real}{\mathbb R}
\newcommand{\Cplx}{\mathbb C}
\newcommand{\Htwo}{\mathbb H^2}

\newcommand{\To}{\longrightarrow}

\def\Set#1{\left\{\,#1\,\right\}}
\def\re{\operatorname{Re}}
\def\im{\operatorname{Im}}
\def\R{\mathbb R}
\begin{document}

\title[Images of harmonic maps]{Images of harmonic maps with symmetry}%
\author{Thomas K. K. Au \& Tom Y. H. Wan}%
\address{}%
\email{}%

\thanks{The second author is partially supported by
   Earmarked Grants of Hong Kong CUHK4291/00P}

\begin{abstract}
We show that under certain symmetry, the images of complete harmonic embeddings
from the complex plane into the hyperbolic plane is completely determined by
the geometric information of the vertical measured foliation and is independent
of the horizontal measured foliation of the corresponding Hopf differentials.
\end{abstract}
\maketitle

\setlength{\baselineskip}{24pt}

In this paper,  we find a new explicit relation between the image of harmonic
embeddings, with certain symmetry, from the complex plane~$\Cplx$ into the
hyperbolic plane~$\Htwo$ and the metric of the associated $\Real$-tree of the
corresponding vertical measured foliation of the Hopf differentials. Unlike in
the case of compact surfaces, holomorphic quadratic differentials cannot be
determined by the vertical measured foliation only. So it is kind of surprising
for us to find that the image set of the corresponding complete harmonic
embedding is completely determined by the vertical measured foliation and is
independent of the geometric information of the horizontal measured foliation.

The symmetry condition that we consider is as follow. We assume that the
harmonic embedding $u$ from $\mathbb{C}$ into $\mathbb{H}^2$ is invariant under
the group $\mathbb{Z}_k$ by rotations and its image is an ideal polygon with
$2k$ vertices for any integer $k\ge 2$. This is the next nontrivial case after
the case of $\mathbb{Z}_{2k}$ symmetry which gives harmonic embeddings with
regular polygonal images. This condition can be regarded as $u$ having half of
the symmetry of a regular polygon.

The symmetry assumption implies that the Hopf differentials are of the form
$[z^{2m}-(a+ib)z^{m-1}]dz^2$ for $a+ib\in \mathbb{C}$. For a generic
holomorphic quadratic differential in this family, the associated $\Real$-tree
has $m+1$ finite edges of equal length given by $\nu=\pi|b|/(2(m+1))$.  We will
show that
\begin{thm*}
Let $u:\mathbb{C}\to\mathbb{H}^2$ be the unique {\em({\em up to
equivalence\/})} complete orientation preserving harmonic embedding associated
to a quadratic differential equivalent to $[z^{2m} -(a+ib)z^{m-1}]dz^2$. Then,
up to isometry, the image $u(\mathbb{C})$ is the interior of the ideal polygon
with vertices given by $\{ 1, e^{i\alpha}, \omega , \omega e^{i\alpha}, \ldots,
\omega^m , \omega^m e^{i\alpha}\}$ in the unit disc model of $\mathbb{H}^2$,
where $\omega=e^{2\pi i/(m+1)}$,
$$
\alpha=\alpha_m(\nu)=2\tan^{-1}\left( \frac{\sin(\pi/(m+1))}{\cos(\pi/(m+1))
+e^{2\nu}} \right),
$$
and $\nu=\pi|b|/(2(m+1))$ is the common length of the finite edges of the
$\mathbb{R}$-tree associated to the quadratic differential given by Lemma~{\em
\ref{lem-tree}}.
\end{thm*}
In this paper, a harmonic embedding $u$ is called {\em complete\/} if its
$\partial$-energy metric $\|\partial u\|^2|dz|^2$ is a complete metric on
$\mathbb{C}$, where $z$ is the standard complex coordinate on $\mathbb{C}$.

The result is related to the work of Shi and Tam \cite{S-T}. The facts that
complete harmonic embeddings from $\mathbb{C}$ to $\mathbb{H}^2$ are
parametrized by Hopf differentials \cite{Tam-Wan, Wan-Au} and the images are
determined by the asymptotic behaviors of the harmonic embeddings \cite{A-T-W,
Han, HTTW}, suggest the following problem as a step toward Schoen's conjecture
\cite{Sch} on the nonexistence of harmonic diffeomorphism from the complex
plane to the hyperbolic plane: Suppose that $u$ is a complete orientation
preserving harmonic embedding with polynomial Hopf differential $P(z)dz^2$, is
it possible to find explicit relation between the coefficients of $P(z)$ and
the vertices of $\overline{u(\mathbb{C})}$? For this problem, they showed that,
up to isometry, the image of a complete orientation preserving harmonic
embedding from the complex plane into the hyperbolic plane is a regular ideal
polygon if its Hopf differential is given by $(z^{2m}-az^{m-1})dz^2$ for some
{\em real} number $a$. This is the first nontrivial example of a family of
harmonic maps (for fixed $m$) with identical images.

It is obvious that our result is a generalization of that of Shi-Tam. However,
the method is quite different. In \cite{S-T}, the authors studied the
asymptotic behavior of the image of the harmonic maps along euclidean rays to
infinity. Our approach adopts more geometric properties of the Hopf
differential, especially those related to the metric information of the
$\mathbb{R}$-tree associated to the vertical measured foliation of the Hopf
differential. The relationship between the asymptotic behavior of harmonic maps
and the associated $\mathbb{R}$-trees has been studied by Minsky \cite{Minsky}
and Wolf \cite{Wolf1, Wolf2, Wolf3}, independently. In these works, the
asymptotic behavior of a sequence of harmonic maps on a compact surface with
energy (or the norm of the Hopf differential) going to infinity was studied. In
our case, instead of a sequence of maps, we are interested in the asymptotic
behavior of harmonic maps on a complete noncompact surface as in \cite{HTTW}.
In particular, the asymptotic behavior of the length of the image of a
horizontal trajectory near infinity was studied.  More precisely, it was shown
that the image of a horizontal trajectory is asymptotic to a geodesic; and the
difference between the lengths of this image and the asymptotic geodesic is
actually tending to zero as the $\Phi$-distance is going to infinity.

The arrangement of this paper is as follows. In Section~\ref{sec-H-tree}, we
will give a brief description of harmonic maps, its Hopf differentials and the
geometric information of the $\mathbb{R}$-trees associated to the Hopf
differentials. Then we will study the asymptotic behavior of the image of
horizontal trajectories in Section~\ref{sec-distance}. Finally, we prove our
main result in Section~\ref{sec-image}.

\section{Background}\label{sec-H-tree}

\subsection{Harmonic maps between surfaces}
Let $M$ and $N$ be oriented surfaces with metrics $\rho^2|d z|^2$ and
$\sigma^2|d u|^2$, where $z$ and $u$ are local complex coordinates of $M$ and
$N$, respectively. A $C^2$ map $u$ from $M$ to $N$ is harmonic if and only if
$u$ satisfies
$$
u_{z \bar{z}}+2\left( \log \sigma(u)\right)_{u}
u_{z}u_{\bar{z}}=0.
$$
The Hopf differential $\Phi =\phi(z)d z^2$ of a map $u$ between these surfaces
is defined by $\phi(z)=\sigma^2\left(u(z)\right)u_z(z)\bar{u}_z(z)$. If $u$ is
harmonic, then it is well-known that $\Phi$ is a holomorphic quadratic
differential on $M$.

The $\partial$-{\em energy density} and $\overline{\partial}$-{\em
energy density} of $u$ are defined by
$$
\|\partial u\|^2=\frac{\sigma^2(u)}{\rho^2}|u_z|^2\quad
\mbox{and}\quad \|\overline{\partial}
u\|^2=\frac{\sigma^2(u)}{\rho^2}|u_{\bar{z}}|^2.
$$
In terms of the $\partial$-{\em energy density} and
$\overline{\partial}$-{\em energy density}, the energy density and
Jacobian of $u$ can be written as
$$
e(u)=\|\partial u\|^2+\|\overline{\partial}
u\|^2\quad\mbox{and}\quad J(u)=\|\partial
u\|^2-\|\overline{\partial} u\|^2.
$$

In this paper, we are interested in the case that $M=\Cplx$, $N=\Htwo$, and
that $u:\Cplx \To \Htwo$ is an orientation preserving open harmonic embedding.
In this case, the Jacobian is strictly positive, i.e., $J(u)>0$, and hence
$\|\partial u\|^2>0$. Therefore, one can consider the $\partial$-energy metric
$\|\partial u\|^2|dz|^2$ on the complex plane $\Cplx$. As mentioned in the
introduction, $u$ is called {\em complete\/} if its $\partial$-energy metric
$\|\partial u\|^2|dz|^2$ is a complete metric on $\mathbb{C}$. As the
completeness is only defined for orientation preserving $u$, the term {\em
complete harmonic open embedding\/} implies implicitly that the harmonic
embedding is orientation preserving.

It was shown in \cite{Tam-Wan, Wan-Au} that for each holomorphic quadratic
differential $\Phi=\phi(z)dz^2$ which is not identically zero, there is a
complete harmonic open embedding, unique up to conformal transformations, $u:
\Cplx \To \Htwo$ such that the Hopf differential of $u$ is exactly $\Phi$.

\subsection{Trajectory structures and measured foliations of the Hopf differentials}
Let $\Phi$ be a holomorphic quadratic differential on~$\Cplx$, which is given
in local coordinate~$z$ as $\Phi = \phi(z)\,dz^2$, where $\phi$ is in general a
holomorphic function. For any $z_0\in\Cplx$ with $\phi(z_0)\ne 0$, there is a
choice of a continuous branch of $\sqrt{\phi(z)}$ in a neighborhood~$W$
of~$z_0$. Then for a given base point $z_*\in W$ sufficiently close to $z_0$,
the mapping
$$
\zeta(z) = \int_{z_*}^z \sqrt{\phi(w)}\,dw
$$
is univalent in possibly a smaller neighborhood of~$z_0$ in~$W$. This defines
local charts on $\Set{\phi\ne 0}$ and determines two measured foliations
on~$\Cplx$ with singularities at the zeros of~$\phi$. In particular, the leaves
of them are curves given locally by the sets,
\begin{align*}
\Gamma_\nu &= \Set{z\in W ~;~ \im(\zeta(z)) = \nu}, \qquad
\nu\in\R, \\
\gamma_\mu &= \Set{z\in W ~;~ \re(\zeta(z)) = \mu}, \qquad \mu\in\R.
\end{align*}
Each $\Gamma_\nu$ and $\gamma_\mu$ is called a {\em horizontal trajectory\/}
and {\em vertical trajectory,\/} respectively. The foliations formed by these
curves are called {\em horizontal foliation\/} and {\em vertical foliation\/}
correspondingly. Obviously, the two foliations have orthogonal leaves.
Furthermore, if $z_0\in\Cplx$ is a zero of order~$m$ of~$\phi$, then there
are~$m+2$ horizontal trajectories, as well as vertical ones, limiting to $z_0$.
Therefore, the horizontal and vertical foliations are in fact {\em measured
foliations} with singularities at the zeros of $\Phi$ with natural measures
given by $|d\mbox{Im}\zeta|$ and $|d\mbox{Re}\zeta|$, respectively. We refer
the reader to \cite{Wolf2} for the definition of measured foliation on Riemann
surface in the general situation.

\subsection{The canonical trees associated to the Hopf
differentials} For each $\Phi$, the leaf space of the measured foliation given
by the vertical trajectories has a special 1-dimensional structure called
${\mathbb R}$-tree \cite{Wolf1, Wolf2, Wolf3}.  In this article,  we shall call
it the {\em $\mathbb{R}$-tree associated to\/}~$\Phi$ and denote it by
$T_{\Phi}$, or simply by $T$.

A trajectory that tends to a zero of~$\Phi$ at least in one direction is called
a {\em critical trajectory\/}.  Each connected domain of the complement of all
critical vertical trajectories is sometimes called a {\em vertical domain\/},
which is foliated by non-critical vertical trajectories.

In the particular case of a quadratic differential $\Phi=P(z)dz^2$ for a
polynomial~$P$ of degree~$n$ on~$\mathbb{C}$, according to the global
structural theorem of meromorphic quadratic differentials on compact Riemann
surfaces \cite{Jenkins}, there are generically $2n+1$ vertical domains. Among
these domains, $n+2$~are called {\em end domains\/} and at most $n-1$~are {\em
strip domains\/}. The definition of these two types of domains is given as
follows.

For each vertical domain~$\Omega$, a canonical mapping $z\mapsto \zeta(z)$
sends $\Omega$ one-to-one onto one of the following domains in~$\Cplx$,
\begin{enumerate}
\item a half plane, in such case $\Omega$ is called an {\em end domain};

\item a vertical strip, $\Set{\zeta\in\Cplx ; a < \re(\zeta) < b}$, $a,b\in\R,$
in such case $\Omega$ is called a {\em strip domain}.
\end{enumerate}

Note that the distance on the $\mathbb{R}$-tree~$T$ can be realized in the
following way.  Let $p$, $q$ be two points on~$T$ represented by two leaves
$\gamma_1$ and $\gamma_2$, respectively.  One may construct a sequence of arcs
from $\gamma_1$ to $\gamma_2$ such that each arc lies either in a horizontal
trajectory or a vertical one.  The distance~$d_T(p,q)$ is given by the sum of
the lengths of the horizontal arcs. In particular, if the straight line between
$p$ and $q$ on the $\mathbb{R}$-tree can be represented by a single horizontal
trajectory in an end domain, then $d_T(p,q)$ equals the $\Phi$-length of that
horizontal trajectory.

Consequently, one can see that the associated $\mathbb{R}$-tree $T$ has $n+2$
infinite edges corresponding to the $n+2$ end domains, at most $n-1$ finite
edges corresponding to the strip domains, and with $n$ vertices corresponding
to the zeros counted with multiplicity.

For the special case that $\Phi=(z^{2m} -cz^{m-1})dz^2$ with generic
$c\in\mathbb{C}$, we see that there are $m+1$ non-degenerate vertices
corresponding to the roots of $z^{m+1}-c$ and $m-1$ vertices degenerated to a
single vertex corresponding to $z=0$ if $m\ge 3$. The tree $T$ will further
degenerate if $c$ is real; and for $c=0$, it will completely degenerate to a
single vertex.  More precisely, we have
\begin{lem}\label{lem-tree}
Let $T$ be the $\mathbb{R}$-tree associated to the quadratic differential
$\Phi=[z^{2m}-(a+ib)z^{m-1}]dz^2$, $a+ib \in \mathbb{C}$, and $m\ge 1$.
\begin{enumerate}
\item If $m\ge 2$ and $b\neq 0$, then $T$ has $m+1$ non-degenerate vertices
each incident with two infinite edges\/{\em ;} and all of these vertices are
adjacent to a unique vertex, which is non-degenerate for $m=2$ and degenerate
otherwise, by finite edges of equal length given by $\pi |b|/(2(m+1))$.

\item If $m=1$ and $b\neq 0$, then $T$ has\/ {\em 2} non-degenerate vertices
each incident with two infinite edges\/{\em ;} and they are connected by a
finite edge of length $\pi |b|/2$.

\item If $m\ge 1$ and $b=0$, then $T$ has a unique vertex incident with $2m+2$ infinite edges.
\end{enumerate}
\end{lem}
\begin{rem}
In the case~(2), if we take the mid-point of the 2 vertices as the center of
the $\mathbb{R}$-tree, then the vertices are in distance $\pi |b|/4$ to this
center. This is exactly the same value given by the formula in the case~(1)
with $m=1$. An illustration is given in Figure~1.
\end{rem}

\begin{center}
\parbox{30mm}{\begin{center}\mbox{\epsfxsize=30mm \epsfbox{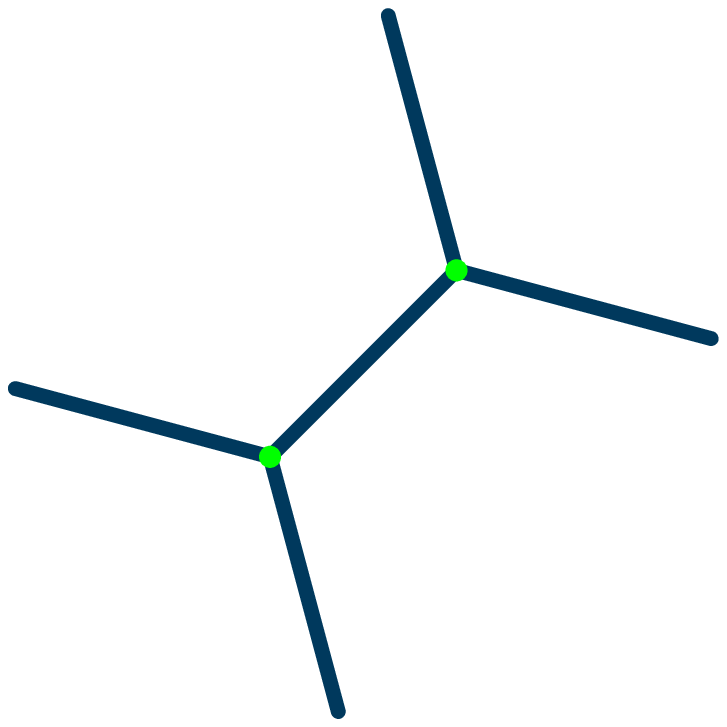}} \\
\end{center}} \hfil
\parbox{30mm}{\begin{center}\mbox{\epsfxsize=30mm \epsfbox{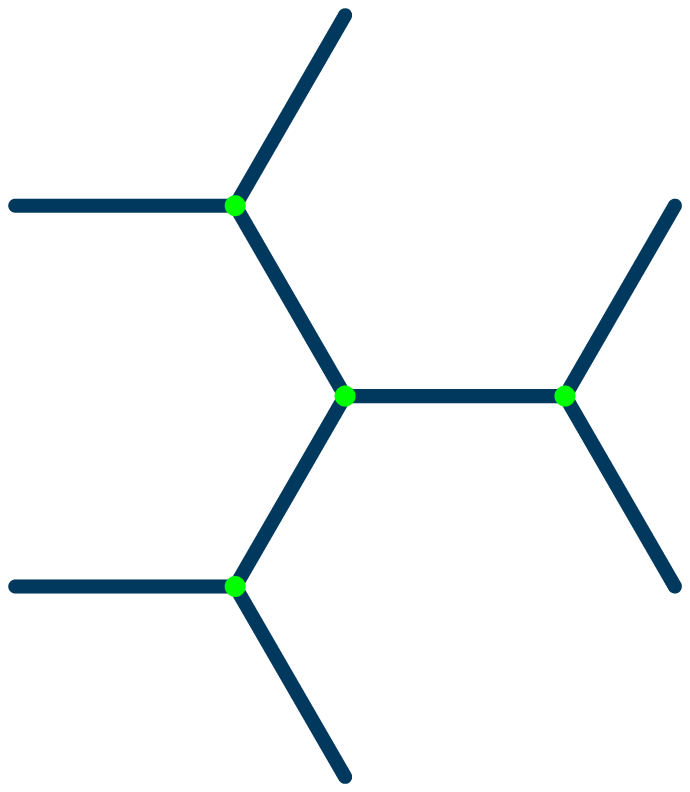}} \\
\end{center}} \hfil
\parbox{30mm}{\begin{center}\mbox{\epsfxsize=30mm \epsfbox{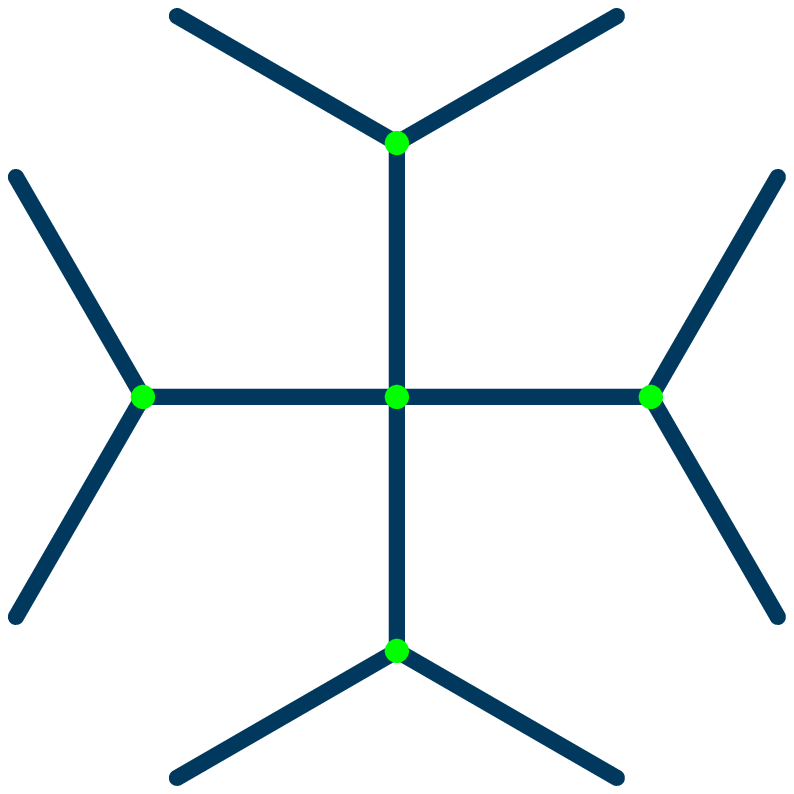}} \\
\end{center}} \\
{\footnotesize Figure 1: $m=1, 2, 3$ from left to right.}
\end{center}

\begin{rem}
If we consider the associated real tree to the leaf space of the measured
foliation given by horizontal trajectories instead of vertical ones, we have
another $\mathbb{R}$-tree associated to $\Phi$. One sees that the same is true
for this tree except that the common length of the finite edges becomes $\pi
|a|/(2(m+1))$.  This tree will degenerate when $c$ is pure imaginary.
\end{rem}
\renewcommand\qedsymbol{}
\begin{proof}
By the argument before the Lemma, $T$ has $2m+2$ infinite edges corresponding
to the $2m+2$ end domains and has $m+1$ non-degenerate vertices. If $m\ge 2$,
then $T$ also has a vertex with multiplicity $m-1$. So we only need to show
that each non-degenerate vertex is incident with two infinite edges and
calculate the distance between the non-degenerate vertices and the vertex with
multiplicity.

For any $m\ge 1$, let $\eta$ be a $(m+1)$-root of $c=a+ib$ and
$\omega$ be a primitive $(m+1)$-root of unity. Then the roots of
$z^{m+1}-c$ are exactly $\{ \eta,\, \omega \eta,\,\ldots,
\omega^{m}\eta\}$. For any fixed $k=0, 1, \ldots, m$, there is a
wedge with vertex at $z=0$, containing the path
$z_k(t)=t\omega^k\eta$, $t\in (0,1)$, but no other zero of $\Phi$.
Therefore, one can find a domain $\Omega_k$ containing the path
$z_k(t)$ which is contained in a strip domain of $\Phi$. Choosing
a branch of the natural parameter
$$
\zeta_k(z)= \int^z\sqrt{z^{m-1}(z^{m+1}-c)}dz
$$
on $\Omega_k$ and taking limits at $t\to 0$ and $t\to 1$, one sees that
$$
\zeta_k(\omega^k\eta)-\zeta_k(0)=\pm ic\int_0^1\sqrt{t^{m-1}(1-t^{m+1})}dt.
$$
This implies that the horizontal $\Phi$-distance between the root
$\omega^k\eta$ and $0$ is given by
$|b|\displaystyle\int_0^1\sqrt{t^{m-1}(1-t^{m+1})}dt=\pi |b|/(2(m+1))$. This
proves that, in the ${\mathbb R}$-tree, the vertices corresponding to the roots
of $z^{m+1}-c$ are adjacent to the vertex corresponding to $z=0$ by a finite
arc of length $\pi |b|/(2(m+1))$ if $m\ge 2$. If $m=1$, then the same
calculation shows that the 2 vertices are connected by a finite arc of length
$2\times \pi |b|/(2(m+1))=\pi |b|/2$.

Therefore, in the case that $m\ge 2$, the vertex corresponding to $z=0$ with
multiplicity $m-1$ already has $m+1$ finite edges incident with the $m+1$
non-degenerate vertices, and hence it is not incident with any other edges of
$T$. Since there are $2m+2$ infinite edges, each of the vertices corresponding
to the roots of $z^{m+1}-c$ must be incident with 2 infinite edges. The case
that $m=1$ is trivial, since $z=0$ is not even a critical point of $\Phi$.
Finally, by counting the multiplicity of the vertices and the number of edges,
we see that there is no other edge of $T$ and the proof is completed.
\end{proof}
\renewcommand\qedsymbol{$\square$}

\section{Distance estimates}\label{sec-distance}

In \cite{HTTW}, it was shown that the image of an horizontal trajectory far
from zeroes is exponentially close to the geodesic connecting the end points in
the $\Phi$-distance of the trajectory. However, this is not enough in our
discussion about the image of the harmonic map. In fact, we need to show that
the difference between the lengths of the image and the geodesic is actually
tending to zero as the $\Phi$-distance is going to infinity. Note that we need
more than just the ratio tending to~$1$ as in Proposition~2.2 of \cite{S-T}.
\begin{lem}\label{lem-dist}
Let $\gamma_R$, $R>0$ be a family of curves in the hyperbolic\/ {\em 2}-space
such that, as $R\to +\infty$, $L(\gamma_R)=O(R)$ and $\|
k_g\|(\gamma_R)=O(e^{-aR})$ for some $a>0$, where $L(\gamma_R)$ is the length
of $\gamma_R$ and $\| k_g\|(\gamma_R)$ is the supremum of the absolute value of
the geodesic curvature of $\gamma_R$. Then the distance $d(R)$ between the end
points of $\gamma_R$ satisfies $d(R)=L(\gamma_R)+o(R)$ as $R\to +\infty$.
\end{lem}
\begin{proof}
For a sufficiently large fixed $R>0$, we work in the Fermi coordinates $(u,v)$
with respect to the geodesic $\gamma_R^{*}$ passing through the end points
$\gamma_R(0)$ and $\gamma_R(l)$, where $l=L(\gamma_R)$ is the length of
$\gamma_R$.  That is, $\gamma_R:[0,l] \to \mathbb{H}^2$ is parametrized by
arc-length, and $\gamma_R^{*}$ is given by $v\equiv 0$.

By Lemma~3.1 in \cite{HTTW},  there exists a constant $C>0$ such that, for
sufficiently small $\epsilon>0$, $\|k_g\|(\gamma_R)<\epsilon$ implies
$d(\gamma_R, \gamma_R^{*})<C\epsilon$. That is,
\begin{equation}\label{eqn-estim1}
|v(\gamma_R(s))|<C\epsilon \quad\mbox{for all } s\in [0,l].
\end{equation}
As in \cite{HTTW}, we have
\begin{equation} \label{eqn-norm1}
{u'}^2\cosh^2v +{v'}^2 \equiv1,
\end{equation}
and
\begin{equation}\label{eqn-geocurv}
k_g^2= \cosh^2v ( u''+2u'v'\tanh v)^2 + ( v''-{u'}^2 \cosh v\sinh v  )^2.
\end{equation}
Let
$$
 h_1 =u'' \cosh v  +2u'v'\sinh v, \qquad
 h_2 =v''-{u'}^2 \cosh v\sinh v.
$$
Then $h_1^2+h_2^2=k_g^2$. On the other hand, differentiation of (\ref{eqn-norm1}) gives $h_1
u'\cosh v + h_2v' =0$, i.e., $(h_1, h_2)$ is orthogonal to $(u'\cosh v, v')$. Therefore, we must
have $(h_1, h_2)=\pm |k_g|(v', -u'\cosh v)$. Consequently, we have
\begin{equation}\label{eqn-estim2}
|h_1|\le \epsilon |v'| \quad\mbox{and}\quad |h_2|\le \epsilon|u'|\cosh v.
\end{equation}
We may assume that $u(0)=0$. Then $d(R)=u(l)$. Applying Poincar\'{e} inequality to $u'(s)-u(l)/l$
and $v'(s)$, we conclude that there is a constant $C_1>0$ such that
\begin{equation}\label{eqn-upoincare}
\int_0^l {u'}^2 \le \frac{d(R)^2}{l} + C_1 l^2 \int_0^l (u'')^2
\end{equation}
and
\begin{equation}\label{eqn-vpoincare}
\int_0^l {v'}^2 \le  C_1 l^2 \int_0^l (v'')^2.
\end{equation}
Therefore, (\ref{eqn-estim1}), (\ref{eqn-norm1}), (\ref{eqn-estim2}) and (\ref{eqn-upoincare})
imply
\begin{equation}\label{eqn-estim3}
\int_0^l {u'}^2 \le \frac{d(R)^2}{l} + C_2 l^2\epsilon^2 \int_0^l {v'}^2
\end{equation}
for some constant $C_2>0$. Similarly, we have from (\ref{eqn-estim1}), (\ref{eqn-norm1}),
(\ref{eqn-estim2}) and (\ref{eqn-vpoincare}) that
\begin{equation}\label{eqn-estim4}
\int_0^l {v'}^2 \le   C_3 l^2\epsilon^2 \int_0^l {u'}^2
\end{equation}
for some constant $C_3>0$. Putting this into (\ref{eqn-estim3}), we have
\begin{equation}\label{eqn-estim5}
\int_0^l {u'}^2 \le \frac{d(R)^2}{l} + C_2C_3 l^4\epsilon^4 \int_0^l {u'}^2.
\end{equation}
By the assumption on the geodesic curvature $\| k_g\|$, we may choose
$\epsilon=O(e^{-aR})$. Then, together with $l=O(R)$, one has $C_2C_3
l^4\epsilon^4 \le\epsilon^2 <1$ for sufficiently large $R$. Hence,
(\ref{eqn-estim5}) gives
\begin{equation}\label{eqn-estim6}
\int_0^l {u'}^2 \le \frac{1}{(1-\epsilon^2)}\frac{d(R)^2}{l} = [1+O(e^{-2aR})]\frac{d(R)^2}{l}.
\end{equation}
On the other hand, from (\ref{eqn-estim1}) and (\ref{eqn-norm1}), we have
$$
l= \int_0^l {u'}^2\cosh ^2 v  +{v'}^2 \le (1+C^2\epsilon^2) \int_0^l {u'}^2 +\int_0^l {v'}^2.
$$
Together with the estimates (\ref{eqn-estim4}) and (\ref{eqn-estim6}), this gives
$$
l\le [1+O(e^{-aR})]\frac{d(R)^2}{l}.
$$
Therefore, $l=O(R)$ implies
$$
d(R)\ge l \left(1+O(e^{-aR})\right)^{-1} = l - l\cdot O(e^{-aR}) \ge l-O(e^{-aR/2}).
$$
As it is trivial that $l\ge d(R)$, we have shown that $d(R)=L(\gamma_R)+o(R)$.
\end{proof}
From Lemma~\ref{lem-dist}, we have the following corollary on the asymptotic
behavior of harmonic maps.
\begin{cor}\label{cor-dist}
Let $\Gamma_R$, $R>0$, be a family of horizontal trajectories of a holomorphic
quadratic differential $\Phi$ with $\Phi$-length equal to $L$. If the
$\Phi$-distance of $\Gamma_R$ to every zero of~$\Phi$ tends to infinity as
$R\to\infty$, then the images $u(\Gamma_R)$ of $\Gamma_R$ under the unique
harmonic embedding $u:\mathbb{C}\to\mathbb{H}^2$ corresponding to $\Phi$
approaches a boundary geodesic arc of length $2L$ of the image set
$u(\mathbb{C})$ as $R\to +\infty$.
\end{cor}
\begin{proof}
Let us write $\gamma_R=u(\Gamma_R)$ and $d_{\phi}(R)$ for the minimal
$\Phi$-distance of $\Gamma_R$ to zeroes of $\Phi$. First of all, the arguments
of Lemma~3.2 and~3.4 of \cite{HTTW} imply that $\gamma_R$ approaches the
boundary geodesic of $u(\mathbb{C})$. So we only need to calculate its length.
Let $e^{2w}$ be the $\partial$-energy density of $u$ with respect to the
$\Phi$-metric in its natural coordinates, i.e., $\Phi=dz^2=(dx+idy)^2$. Then,
by Formula~(3.6) of \cite{HTTW}, we obtain
$$
L(\gamma_R)=\int_{\Gamma_R} \sqrt{e^{2w}+e^{-2w}+2}dx =\int_{\Gamma_R} \left(2+4\sinh
\frac{w}{2}\right) dx.
$$
The exponential decay estimate of \cite{Han} then implies that
$$
L(\gamma_R)=2L + O(e^{-a_1d_{\Phi}(R)})
$$
for some constant $a_1>0$. On the other hand, the estimate as in the proof of
Lemma~3.2 in \cite{HTTW} shows that
$$
\|k_g\|(\gamma_R)=O(e^{-a_2d_{\Phi}(R)})
$$
for some constant $a_2>0$. Therefore, by Lemma~\ref{lem-dist}, we conclude that
the distance between the end points of $\gamma_R$ is equal to $2L+o(R)$ as
$R\to \infty$. Therefore, by letting $R\to\infty$, we have the desired result.
\end{proof}

\section{Image of harmonic maps}\label{sec-image}

In this section, we prove our main result on the explicit determination of the
image of the harmonic embedding with suitable symmetry. We are interested in a
harmonic embedding $u$ from $\mathbb{C}$ into $\mathbb{H}^2$ such that $u$ is
equivariant under the group $\mathbb{Z}_k$ by rotations and its image is an
ideal polygon with $2k$ vertices for any integer $k\ge 2$. In some sense, $u$
has half of the symmetry of a regular polygon. Note that our symmetry
requirement is not just on the image set but on the map $u$.

According to this requirement, the Hopf differentials of these harmonic
embeddings are equivariant under the action $z\mapsto \omega z$ for any
$k^{\text{th}}$-root of unity~$\omega$ and their coefficients are polynomials
of degree $2k-2$. This immediately implies that the Hopf differentials are of
the form $[z^{2m} -(a+ib)z^{m-1}]dz^2$, where $a+ib \in\mathbb{C}$. For these
type of harmonic embeddings, we have the following
\begin{thm}\label{thm-main}
Let $u:\mathbb{C}\to\mathbb{H}^2$ be the unique\/ {\em ({\em up to
equivalence\/})} harmonic embedding associated to a quadratic differential
equivalent to $[z^{2m} -(a+ib)z^{m-1}]dz^2$. Then, up to isometry, the image
$u(\mathbb{C})$ is the interior of the ideal polygon with vertices given by
$$
\left\{ 1, e^{i\alpha}, \omega , \omega e^{i\alpha}, \ldots, \omega^m ,
\omega^m e^{i\alpha}\right\}
$$
in the unit disc model of\/ $\mathbb{H}^2$, where $\omega=e^{2\pi i/(m+1)}$,
$$
\alpha=\alpha_m(\nu)=2\tan^{-1}\left( \frac{\sin(\pi/(m+1))}{\cos(\pi/(m+1))
+e^{2\nu}} \right),
$$
and $\nu=\pi|b|/(2(m+1))$ is the common length of the finite edges of the
$\mathbb{R}$-tree associated to the quadratic differential given by Lemma~{\em
\ref{lem-tree}}.
\end{thm}
\begin{proof}
As a harmonic map from a surface is invariant under conformal change of metrics
on the surface, we may assume that the Hopf differential of $u$ is in fact
given by $\Phi=[z^{2m} -(a+ib)z^{m-1}]dz^2$. Then, by the symmetry of the
quadratic differential $\Phi$ and the uniqueness property of the corresponding
complete orientation preserving harmonic embedding, after a composition with an
isometry on $\mathbb{H}^2$, the harmonic embedding satisfies $u(0)=0$ and the
image $u(\mathbb{C})$ is an ideal polygon with vertices given by
$$
\{ 1, e^{i\alpha}, \omega , \omega e^{i\alpha}, \ldots, \omega^m , \omega^m
e^{i\alpha}\}
$$
for some $\alpha\in (0, 2\pi/(m+1))$. What we need to do is to determine
$\alpha$. We also note that, by the rotation of an angle~$-\alpha$, this
polygon is equivalent to $\{ 1, e^{i\beta}, \omega , \omega e^{i\beta}, \ldots,
\omega^m , \omega^m e^{i\beta}\}$ with $\beta=2\pi/(m+1)-\alpha$. Therefore, we
may assume that $\alpha\in(0,\pi/(m+1))$.

Let $0\in (T, d_T)$ be the vertex on the associated $\mathbb{R}$-tree not
incident with any infinite edge for $m\ge 2$ or the mid-point of the unique
pair of vertices for $m=1$ as described in Lemma~\ref{lem-tree}. For
sufficiently large $L>0$, the set $\{ q\in T \,;\, d_T(q, 0)=L \}$ has exactly
$2m+2$ points $\{ q_0, \ldots, q_{2m+1}\}\subset T$ such that each infinite
edge contains exactly one $q_i$. As the tree $T$ is coming from the
trajectories structure of $\Phi$ on the plane, there is a natural induced
cyclic order of the set of infinite edges. Assume that $q_i$ are labelled in
the same cyclic order, for $i\in{\mathbb Z}_{2m+2}$.  Then for each pair of
consecutive points $q_i, q_{i+1}$, we can find $z_i$ and $z'_{i+1}$ in
$\mathbb{C}$ both contained in a common horizontal trajectory $\Gamma_i$ of
distance $R$ to zeroes in an end domain, denoted by $E_i$, of $\Phi$ such that
each $z_i$ and $z'_{i+1}$ belongs to the vertical trajectories representing
$q_i$ and $q_{i+1}$, respectively. Note that from our choice, $z_i$ and $z_i'$
belong to the same vertical trajectory $\gamma_i$ representing $q_i$. An
illustration is given in Figure~2.

\begin{center}
\mbox{\epsfysize=60mm \epsfbox{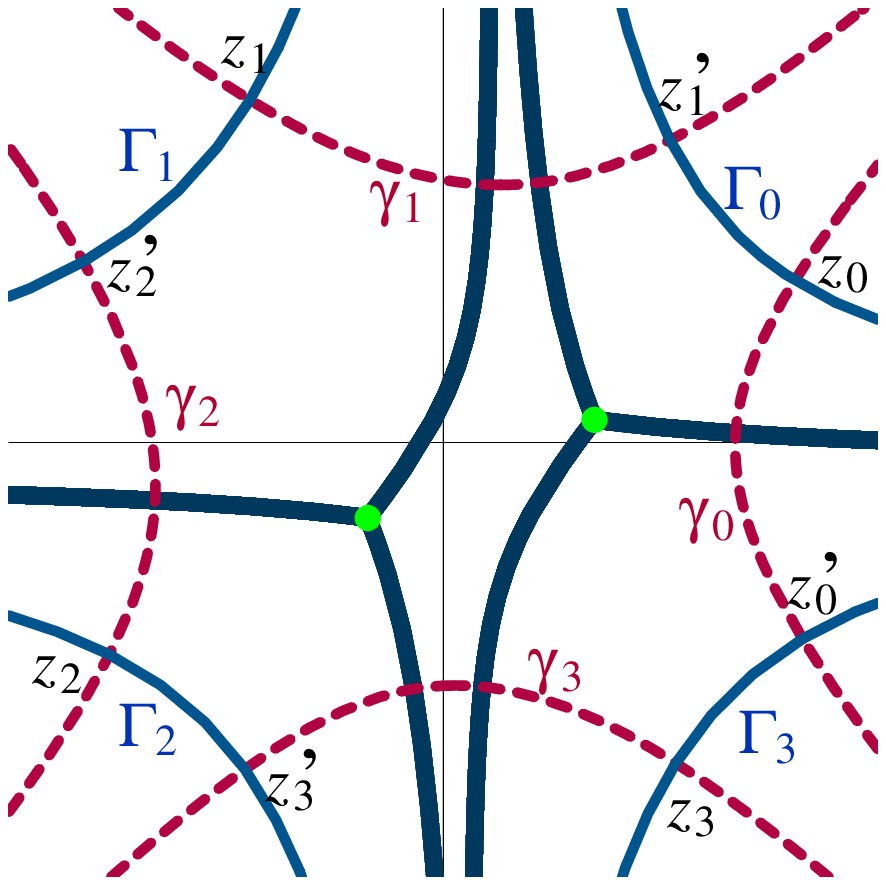}} \\
{\footnotesize Figure 2.}
\end{center}

Up to isometry, we may assume that the image curve of the vertical trajectory
representing $q_0$ approaches the ideal boundary to the point $1$ in the unit
disc model of $\mathbb{H}^2$ as $L\to+\infty$. Correspondingly, the image
points $u(z_0)$ and $u(z_{0}')$ of $z_0$ both tend to $1$. Then, by the
symmetry of $u$ and our assumption, for each $k=0,\ldots, m$, the image curve
of the vertical trajectory representing $q_{2k}$ approaches the ideal boundary
point $\omega^k=e^{2\pi ki/(m+1)}$, and the image curve of the vertical
trajectory representing $q_{2k+1}$ approaches the ideal boundary point
$\omega^ke^{i\alpha}$, respectively, in the unit disc model of $\mathbb{H}^2$.
In accordance with this, the image points $u(z_{2k})$ and $u(z_{2k}')$ tend to
$\omega^k$, while $u(z_{2k+1})$ and $u(z_{2k+1}')$ tend to
$\omega^ke^{i\alpha}$.

To determine the $\Phi$-length of each $\Gamma_i$, we observe from
Lemma~\ref{lem-tree}, which concerns the tree structure of $\mathbb{R}$-tree
$T$ associated to $\Phi$, that $d_{T}(q_i, q_{i+1})=2L$ or $2(L-\nu)$, with the
value taken alternatingly in $i$, where $\nu=\pi |b|/(2(m+1))$ is the common
length of those finite edges of $T$. We first assume that
$$
d_T(q_{2k}, q_{2k+1})= 2(L-\nu)\quad\mbox{and}\quad d_T(q_{2k+1},q_{2k+2})=2L.
$$
An illustration is given in Figure~3.
\begin{center}
\mbox{\epsfysize=52mm \epsfbox{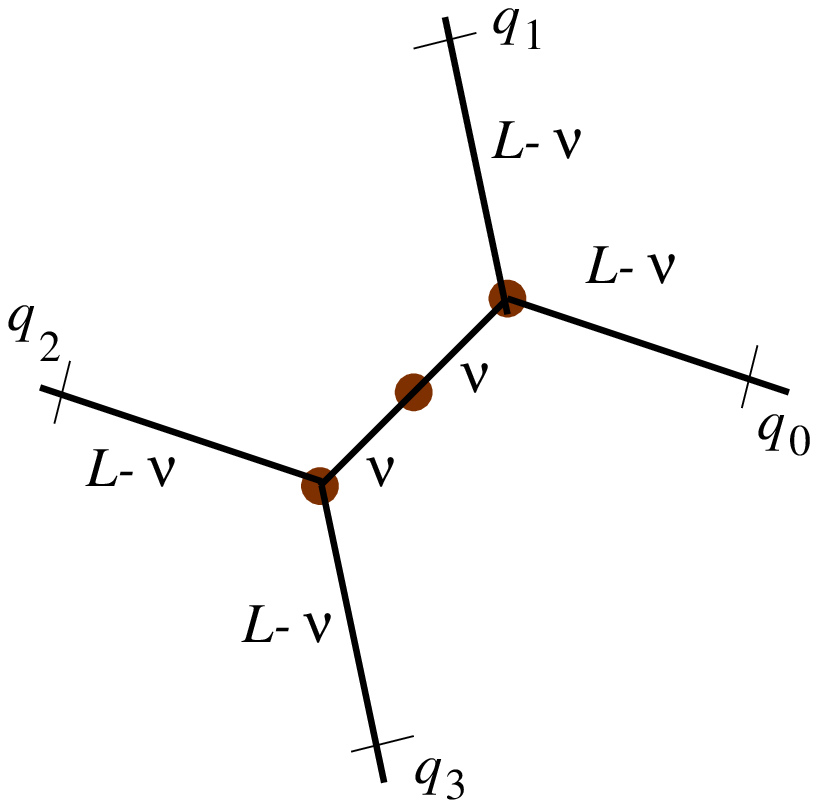}} \\
{\footnotesize Figure 3.}
\end{center}
As $\Gamma_i$ is a horizontal trajectory with end points representing $q_i$ and
$q_{i+1}$, the $\Phi$-length of $\Gamma_i$ is exactly equal to
$d_T(q_i,q_{i+1})$. Therefore, for $k=0,\ldots, m$,
$$
L_{\Phi}(\Gamma_{2k})=2(L-\nu)\quad\mbox{and}\quad L_{\Phi}(\Gamma_{2k+1})=2L.
$$

On the other hand, the vertical trajectory $\gamma_i$ representing $q_i$ is mapped to a curve of
finite length in $\mathbb{H}^2$. Indeed, using natural coordinates of $\Phi$ in an end domain
containing $\gamma_i$ with respect to the vertical trajectories system, the length of the image
curve is given by
$$
l_i=L_{\mathbb{H}^2}(u(\gamma_i))=\int_{-\infty}^{+\infty}\sqrt{e^{2w}+e^{-2w}-2}dy=
\int_{-\infty}^{+\infty}2\sinh w dy,
$$
where $w$ as in the proof of Corollary~\ref{cor-dist}. As $\gamma_i$ is at
least a $\Phi$-distance of $L-\nu$ away from zeroes, the exponential decay
estimate of $w$ implies that for some $y_0$ and $a>0$,
\begin{eqnarray*}
l_i&\le& C \left[ \int_{-y_0}^{y_0} e^{-aL}dy + \left(\int_{-\infty}^{-y_0}  +
\int_{y_0}^{+\infty}\right)
e^{-a(L+|y|)/\sqrt{2}} dy \right] \\
& \le & O(e^{-aL/\sqrt{2}}).
\end{eqnarray*}
Therefore, $l_i$ are finite and tends to zero as $L\to+\infty$.

Let $\zeta_{2k}$ be a point on the intersection of $u(\gamma_{2k})$ and the ray
from $0$ to $\omega^k$ in the Poincar\'{e} disc. Similarly, let $\zeta_{2k+1}$
be a point on the intersection of $u(\gamma_{2k+1})$ and the ray from $0$ to
$\omega^ke^{i\alpha}$. Now consider the polygon in $\mathbb{H}^2$ with vertices
$0, \zeta_0, \zeta_1$, and $\zeta_2$. Note that, since $l_i\to 0$ as
$L\to+\infty$, the distance between $\zeta_i$ and $u(z_i)$ or $u(z_i')$ also
tends to zero. An illustration is given in Figure~4.
\begin{center}
\mbox{\epsfysize=50mm \epsfbox{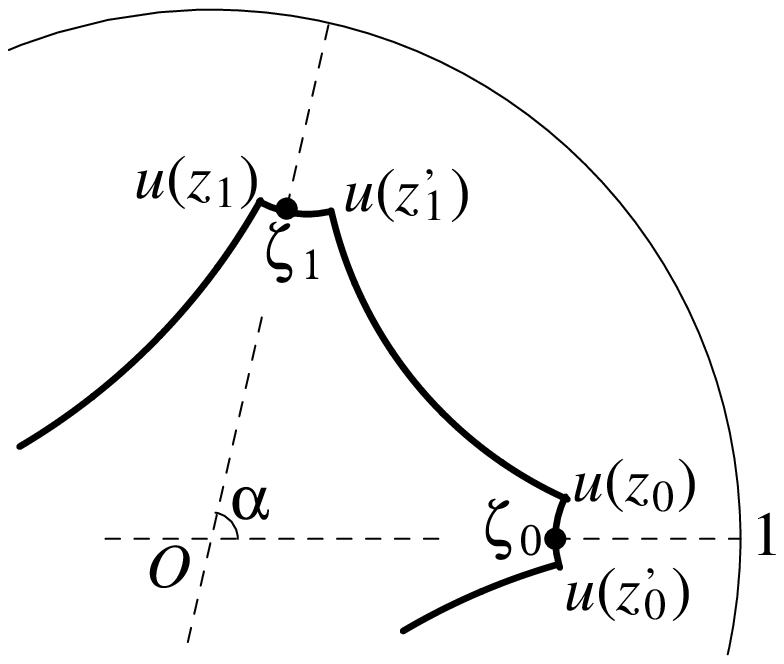}} \\
{\footnotesize Figure 4.}
\end{center}
On the other hand, by letting $R\to+\infty$, Corollary~\ref{cor-dist} implies
that the image $u(\Gamma_i)$ of the horizontal arc $\Gamma_i$ connecting $z_i$
and $z_i'$ approaches a boundary geodesic arc of length $2(L-\nu)$ and $2L$,
alternatingly in $i$. All together, we conclude that, as $L\to+\infty$,
$$
d_{\mathbb{H}^2}(\zeta_0, \zeta_1)=4(L-\nu)+o(L)\quad \mbox{and}\quad d_{\mathbb{H}^2}(\zeta_1,
\zeta_2)=4L+o(L).
$$
Let $x_1=x_1(L)=d_{\mathbb{H}^2}(\zeta_0, 0)$ and
$x_2=x_2(L)=d_{\mathbb{H}^2}(\zeta_1, 0)$. Then, by symmetry,
$d_{\mathbb{H}^2}(\zeta_2, 0)$ is also equal to $x_1$. Hence from the cosine
rule of the hyperbolic plane, we have
$$
\cosh (4(L-\nu)+o(L))=\cosh x_1\cosh x_2 - \sinh x_1\sinh x_2 \cos\alpha
$$
and
$$
\cosh (4L+o(L))=\cosh x_1\cosh x_2 - \sinh x_1\sinh x_2
\cos\left(2\pi/(m+1)-\alpha\right).
$$
It is easy to see from these identities that $\lim_{L\to+\infty}(e^{-4L}\sinh
x_1\sinh x_2)$ exists and is non-zero. Let us denote
$$
A=[{4\lim_{L\to+\infty}(e^{-4L}\sinh x_1\sinh x_2})]^{-1}.
$$
Then multiplying by $\left(\sinh x_1\sinh x_2\right)^{-1}$ to the above
equations and letting $L\to+\infty$, one concludes that
$$
\sqrt{A}e^{-2\nu}=\sin \frac{\alpha}{2}\quad\mbox{and}\quad \sqrt{A}= \sin
\left(\frac{\pi}{m+1}-\frac{\alpha}{2}\right).
$$
It is easy to solve the above and obtain
$$
\tan\frac{\alpha}{2}=\frac{\sin(\pi/(m+1))}{\cos(\pi/(m+1)) +e^{2\nu}}\,,
$$
which is the desired result.

In the case that
$$ d_T(q_{2k}, q_{2k+1})= 2L\quad\mbox{and}\quad d_T(q_{2k+1},q_{2k+2})=2(L-\nu),
$$
the same calculation shows that the angle is given by
$$
\tan\frac{\alpha}{2}=\frac{\sin(\pi/(m+1))}{\cos(\pi/(m+1)) +e^{-2\nu}}.
$$
The angle obtained in this formula belongs to $[\pi/(m+1), 2\pi/(m+1))$, which
is equivalent to the one in previous formula by the transformation $\alpha
\mapsto 2\pi/(m+1)-\alpha$.
\end{proof}

In \cite{S-T}, a harmonic map was constructed with image equal to a regular
ideal polygon of 4 vertices and Hopf differential is given by $(z^2+ib)dz^2$
for some real number $b\in\mathbb{R}$. From our theorem, one in fact has
\begin{cor}
The harmonic diffeomorphism constructed in Proposition~{\em 1.6} of \cite{S-T}
is a unique, up to equivalence, complete orientation preserving harmonic
embedding with Hopf differential $z^2dz^2$.
\end{cor}
\begin{proof}
When $m=1$, the theorem implies that the image of the harmonic map is
equivalent to  $\{ 1, e^{i\alpha}, -1, -e^{i\alpha}\}$ with
$\alpha=2\tan^{-1}(e^{-2\nu})$. So $\alpha=\pi/2$ if and only if $\nu=0$. Since
$\nu=2\pi |b|/(m+1)$, we conclude that $b=0$.
\end{proof}

Finally, let us finish the paper by a couple of remarks.
\begin{rem}
The fact that the image ideal polygon depends only on $|b|$ but not $b$ can be
easily seen from the fact that $[z^{2m}-(a+ib)z^{m-1}]dz^2$ is equivalent to
$[z^{2m}+(a+ib)z^{m-1}]dz^2$.
\end{rem}
\begin{rem}
Let $P_\alpha$ be the equivalence class of the ideal polygon
$$
\left\{ 1, e^{i\alpha}, \omega, \omega e^{i\alpha}, \ldots, \omega^m, \omega^m e^{i\alpha}
\right\}.
$$
Then the mapping $ \alpha \mapsto P_\alpha $ from $\left( 0, 2\pi/(m+1)
\right)$ to all such equivalence classes is two-to-one except at $\alpha =
\pi/(m+1)$, which maps to the regular ideal polygon. Therefore from the proof,
one may define, for each $\Phi=[z^{2m}-(a+ib)z^{m-1}]dz^2$, the angle function
by
$$
\alpha(b)=2\tan^{-1}\left( \frac{\sin(\pi/(m+1))}{\cos(\pi/(m+1))
+e^{b\pi/(m+1)}} \right) ,
$$
that is, by the same formula without taking absolute value of $b$ as in
$\nu=\pi|b|/(2(m+1))$. This angle function~$\alpha$ is a bijection from
${\mathbb R}$ to $\left( 0, 2\pi/(m+1) \right)$.  Thus, the mapping $b \mapsto
P_{\alpha(b)}$ behaves similarly. This gives a 2-fold covering except $b=0$ for
each fixed $a$ and is consistent with the previous remark.
\end{rem}

Note that if we let the angle $\alpha$ run through the whole interval
$(0,2\pi/(m+1))$, the ideal polygon $\{ 1, e^{i\alpha}, \omega , \omega
e^{i\alpha}, \ldots, \omega^m , \omega^m e^{i\alpha}\}$ runs over the set of
all possible equivalence classes of polygons except the regular polygon twice
and once at the regular idea polygon. Then for each fixed $a$, $\alpha :\,
\mathbb{R} \to (0, 2\pi/(m+1))$ is a bijection and the corresponding ideal
polygon with vertices $\{ 1, e^{i\alpha}, \omega , \omega e^{i\alpha}, \ldots,
\omega^m , \omega^m e^{i\alpha}\}$ runs through the set of all possible
equivalent classes of ideal polygons
 except the regular idea polygon twice and once at the regular idea polygon as $b$ run through $\mathbb{R}$ once.


\vskip 2ex

\begin{flushleft}\setlength\baselineskip{12pt}
{\footnotesize {\sc Department of Mathematics \\
The Chinese University of Hong Kong \\
Shatin, Hong Kong \\
China} \\
{\em E-mail addresses\/}: {\tt tomwan@math.cuhk.edu.hk, thomasau@cuhk.edu.hk}}
\end{flushleft}
\end{document}